\documentclass[10pt]{amsart}

\usepackage{amsfonts,verbatim,amsmath,amsthm,amssymb,amscd,enumerate,eucal,url,stmaryrd}

\usepackage[hypertex]{hyperref}      

\numberwithin{equation}{section}
\usepackage{graphicx}
\newtheorem{thrm}{Theorem}[section]
\newtheorem{lemma}[thrm]{Lemma}
\newtheorem{prop}[thrm]{Proposition}
\newtheorem{cor}[thrm]{Corollary}
\newtheorem{dfn}[thrm]{Definition}

\newtheorem{conv}[thrm]{Convention}

\newcommand{\bG}{\boldsymbol {G}}

\newcommand{\Cn}{\mathbb{C}^n}
\newcommand{\Om}{\Omega}

\newcommand{\QH}{\boldsymbol {G\,(\mathbb{P})}}

\newcommand{\lap}{\mathcal{L}}

\newcommand{\st}{\text{*}}

\newcommand{\abs}[1]{\lvert #1 \rvert}

\newcommand{\Ab}{\Bar{A}}
\newcommand{\nz}{\lvert z \rvert}
\newcommand{\nw}{(|u|^2-|v|^2)}

\newcommand{\nB}{\lvert B \rvert}
\newcommand{\vu}{\mathbf{u}}
\newcommand{\vv}{\mathbf{v}}
\newcommand{\vx}{\mathbf{x}}
\newcommand{\vy}{\mathbf{y}}
\newcommand{\vw}{\mathbf{w}}

\newcommand{\dxj}{\frac {\partial {}} {\partial x_{j}} }
\newcommand{\dxmj}{\frac{\partial {}}{\partial x_{m+j}}}
\newcommand{\dyj}{\frac{\partial {}}{\partial y_{j}}}
\newcommand{\dymj}{\frac{\partial {}}{\partial x_{m+j}}}
\newcommand{\ddt}{\frac{\partial {}}{\partial t}}

\newcommand{\Cr}{\nabla^{cr}}
\newcommand{\Rcr}{R^{cr}}
\newcommand{\GP}{G\,(\mathbb{P})}

\begin{document}

\begin{abstract}
A  curvature-type tensor invariant called para contact (pc) conformal
curvature is defined on a paracontact manifold.
It is shown that a paracontact manifold is locally paracontact
conformal to the hyperbolic Heisenberg group or to a hyperquadric of
neutral signature if and only if  the pc conformal curvature vanishes. In the three
dimensional case the corresponding result is achieved through employing a
certain symmetric (0,2) tensor.
 The well known result of Cartan-Chern-Moser giving necessary and
sufficient condition a CR-structure to be CR equivalent to a hyperquadric
in $\mathbb{C}^{n+1}$ is presented in-line with the paracontact case.
An explicit formula for the regular part of a solution to the
sub-ultrahyperbolic Yamabe equation on the hyperbolic Heisenberg
group is shown.
\end{abstract}

\title[Conformal paracontact curvature] 
{Conformal paracontact curvature\\ and the local flatness theorem}
\date{\today}
\thanks{This project has been funded in part by the National Academy
of Sciences under the [Collaboration in Basic Science and Engineering Program 1 Twinning
Program] supported by Contract No. INT-0002341 from the National Science Foundation. The
contents of this publication do not necessarily reflect the views or policies of the
National Academy of Sciences or the National Science Foundation, nor does mention of
trade names, commercial products or organizations imply endorsement by the National
Academy of Sciences or the National Science Foundation.}
 \author{Stefan Ivanov}
\address[Stefan Ivanov]{University of Sofia, Faculty of Mathematics and Informatics,
blvd. James Bourchier 5, 1164, Sofia, Bulgaria}


\email{ivanovsp@fmi.uni-sofia.bg}

\author{Dimiter Vassilev}
\address[Dimiter Vassilev]{
Department of Mathematics and Statistics \\ University of New
Mexico\\ Albuquerque, New Mexico 87131}


\email{vassilev@math.unm.edu} 


\author{Simeon Zamkovoy}
\address[Simeon Zamkovoy]{University of Sofia, Faculty of Mathematics and Informatics,
blvd. James Bourchier 5, 1164, Sofia, Bulgaria}

\email{zamkovoy@fmi.uni-sofia.bg}

\maketitle

\setcounter{tocdepth}{2}


\section{Introduction}

A paracontact structure on a real (2n+1)-dimensional manifold $M$
is a codimension one distribution $\mathbb H$ and a paracomplex
structure $I$ on $\mathbb H$, i.e. $I^2=id$ and the $\pm$
eigen-distributions  $H^{\pm}$ have equal dimension.  Locally,
$\mathbb H$ is the kernel of a 1-form $\eta, \, \mathbb
H=Ker\,\eta$. A paracontact Hermitian structure is a paracontact
structure with the additional assumption that $\eta$ is a para
Hermitian contact form in the sense that we have a non-degenerate
pseudo-Riemannian metric $g$, which is defined on $\mathbb H$, and
compatible with $\eta$ and $I,\quad d\eta(X,Y)=2g(IX,Y),\quad
g(IX,IY)=-g(X,Y),\quad X,Y\in \mathbb H
$. 
The signature of $g$ on $\mathbb H$ is necessarily of (signature)
type (n,n). A para contact structure is said to be integrable if
the para complex structure $I$ on $H$ is formally integrable,
i.e., $[ H^{\pm},H^{\pm}]\subset H^{\pm}.$ A para-contact manifold
with an integrable para-contact structure is called a \emph{para
CR-manifold}

The 1-form $\eta$ is determined up to a conformal factor and hence
$\mathbb H$ become equipped with a conformal class [g] of neutral
Riemannian metrics of signature (n,n). Transformations preserving a
given para contact hermitian structure $\eta$, i.e.
$\bar\eta=\mu\eta$ for a non-vanishing smooth function $\mu$ are
called \emph{para contact conformal (pc conformal for short)
transformations} \cite{Zam}.

A basic example is provided by a para-Sasakian manifold, which can be
defined as a $(2n+1)$-dimensional Riemannian manifold whose metric cone is
a para-K\"ahler manifold \cite{ACGL}. It was shown in \cite{Zam} that the
torsion endomorphism of the canonical connection is the obstruction for a
given integrable para contact hermitian structure to be locally
para-Sasakian, up to a multiplication with a constant factor.

Any non-degenerate hypersurface in $(R^{2n+2},\mathbb I,g)$
considered with the standard flat parahermitian structure inherits
an integrable para-contact hermitian structure. We consider the
(2n+1)-dimensional Heisenberg group with a left-invariant
para-contact hermitian structure $\eta$ and call it
 \emph{hyperbolic Heisenberg group}, denoted by $(\QH,\eta)$. We show that in dimension greater
 than three the hyperbolic Heisenberg group is the unique example of
 an integrable
para-contact hermitian structure with flat canonical connection.
In the three dimensional case the same statement holds under the
additional assumption of vanishing of the torsion tensor. {As a
manifold $\QH \ =\mathbb R^{2n}\times \mathbb R$ with the group
law given by
$
( p'',t'')\ =\ (p',t')\circ(p, t)\ =\ (p'\ +\ p, t' +\ t\ -
\sum_{k=1}^n(u'_kv_k-v'_ku_k),
$
where $p',p\in\mathbb R^{2^n}$ with the standard
coordinates $(u_1,v_1,\dots,u_n,v_n)$ and $t', t\in \mathbb R$.
Define the 'standard' para-contact structure by the left-invariant
para-contact form
$$\tilde\Theta\ =\ -\frac 12\ dt\ - \sum_{k=1}^n(u_kdv_k-v_kdu_k).$$

In this paper we find a tensor invariant characterizing locally the
integrable para-contact hermitian structures which are para-contact
conformally equivalent to the flat structure on $\QH$.  To any integrable
para-contact hermitian structure we associate a curvature-type tensor
$W^{pc}$ defined in terms of the curvature and torsion of the canonical
connection by \eqref{pinv}, whose form is similar to the Weyl conformal
curvature in Riemannian geometry (see e.g. \cite{Eis}) and to the
Chern-Moser curvature in CR geometry \cite{ChM}. We call $W^{pc}$
\emph{para-contact conformal curvature or pc conformal curvature}.
When $M$ is three dimensional, we define in \eqref{treeF} a
symmetric (0,2) tensor $F$ on $\mathbb H$, which plays a role
similar to the Schouten tensor in a 3-dimensional locally
conformally flat Riemannian manifold.

The main purpose of this article is to prove the following two
results.
\begin{thrm}\label{main1}
The pc conformal curvature $W^{pc}$ of an integrable para-contact
hermitian  manifold is invariant under para-contact conformal
transformations.
\end{thrm}
\begin{thrm}\label{main2} Let $(M,\eta)$ be a $2n+1$ dimensional integrable para-contact hermitian  manifold.
\begin{enumerate}
\item[i)] If $n>1$ then $(M,\eta)$ is locally para-contact
conformal to the standard flat para-contact hermitian structure on
the hyperbolic Heisenberg group $\QH$ if and only if the
para-contact conformal curvature vanishes, $W^{pc}=0$. \item[ii)]
If $n=1$ then $W^{pc}$ vanishes identically and $(M,\eta)$ is
locally para-contact conformal to the standard flat para-contact
hermitian structure on the 3-dimensional hyperbolic Heisenberg
group $\QH$ if and only if the  symmetric tensor $F$ vanishes,
$F=0$.
\end{enumerate}
\end{thrm}
We define a Cayley transform which establishes a conformal para-contact
equivalence between the standard para-Sasaki structure on the hyperboloid,
cf. \ref{standard hyperboloid},
\begin{equation}\label{e:ex_hyperboloid}
HS^{2n+1}=\{(x_1,y_1,\dots,x_{n+1},y_{n+1}):
x^2_1+\dots+x^2_{n+1}-y^2_1-\dots-y^2_{n+1}=1\} \subset\mathbb
R^{2n+2}\end{equation} and the standard para-contact hermitian  structure
on $\QH$. As a consequence of Theorem~\ref{main2} and the fact that the
Cayley transform is a para-contact conformal equivalence between the
 hyperboloid  and the group $\QH$, we obtain
\begin{cor}
Let $(M,\eta)$ be a $2n+1$ dimensional integrable para-contact hermitian
manifold. $(M,\eta)$ is locally para-contact conformal to the hyperboloid
$HS^{2n+1}$ if and only if conditions i) or ii) of Theorem \ref{main2}
hold.
\end{cor}
Our investigations are  close to the classical approach used by H. Weyl
(see e.g. \cite{Eis}) and follow the steps of \cite{IV}, compare with
\cite{ChM} where the Cartan method of equivalence is used.
Recall,
that in the CR case  the vanishing of the Chern-Moser tensor is a
necessary and sufficient condition for a non-degenerate CR manifold M of
dimension $2n+1$, $n>1$, to be locally equivalent to a real hyperquadric
in $\mathbb{C}^{n+1}$ of the same signature as M. When $M$ is three
dimensional, the same conclusion can be reached using the Cartan invariant
\cite{Car}.
Both results can be obtained following the steps of the proof of
Theorem~\ref{main2}. In particular, we express the flatness
condition for an abstract three-dimensional pseudohermitian
structure in terms of the covariant derivatives of the
pseudohermitian scalar and torsion of the Webster connection.

{ Let us note that, as observed, the hyperboloid $HS^{2n+1}$ is
always para-contact conformally flat, while the hyperboloid
$HS^{4n+1}$ considered as an embedded CR submanifold of
$\mathbb{C}^{2n+1}$ is a degenerate CR manifold.}

In the last section we consider the CR-Yamabe equation  on a CR
manifold of neutral signature. This leads to the non-linear sub
ultra-hyperbolic equation \eqref{e:Yamabe}, which coincides with the
Yamabe equation for the considered para CR manifolds. Using this
relation we show an explicit formula for the regular part of
solutions to the Yamabe equation.

The paper uses a Webster-like connection, the canonical connection
considered in \cite{Zam} and the properties of its torsion and
curvature described in Section~\ref{bib}.

\begin{conv}\label{conven}
 In the first five sections of the  paper we use
\begin{enumerate}[a)]
\item $X,Y,Z...$ denote horizontal vector fields,
i.e. $X,Y,Z...\in\mathbb H$ \item
$\{e_1,\dots,e_{n},Ie_1,\dots,Ie_{n}\}$ denotes an adapted
orthonormal basis of the horizontal space $\mathbb H$. \item The
summation convention over repeated vectors from the basis
$\{e_1,\dots,e_{2n}\}$ will be used. For example, for a
(0,4)-tensor $P$ we have
\[
P(e_b,X,Y,e_b)=\sum_{b=1}^{n}g(e_b,e_b)P(e_b,X,Y,e_b)+\sum_{b=1}^{n}g(Ie_b,Ie_b)P(Ie_b,X,Y,Ie_b).
\]
\end{enumerate}
\end{conv}

\textbf{Acknowledgements} The research was done during the visit
of S.Ivanov and D. Vassilev in the Max-Plank-Institut f\"ur
Mathematics, Bonn. They thank MPIM, Bonn for providing the support
and an excellent research environment. S.I. is a Senior Associate
to the Abdus Salam ICTP. S.I. and S. Zamkovoy are 
partially supported by the Contract 154/2008 with the University of Sofia "St. Kl. Ohridski".
The authors would like to thank the referee for his useful
suggestions on the presentation of the paper and especially for
his comments concerning the 3-dimensional case { leading to a
clarification of the integrability condition in dimension 3, which
resulted in the current explicit formula.}

\section{Integrable para-contact manifolds}
A para-contact manifold $(M^{2n+1},\eta,I,g)$ is a
(2n+1)-dimensional smooth manifold equipped with a codimension one
distribution $\mathbb H$, locally given as the kernel of a 1-form
$\eta, \, \mathbb H=Ker\,\eta$ and a paracomplex structure $I$ on
$\mathbb H$. Recall that a paracomplex structure is an endomorphism
$I$ satisfying $I^2=id$ and the $\pm$ eigen-distributions  have
equal dimension. If in addition there exists  a  pseudo-Riemannian
metric $g$ defined on $\mathbb H$ compatible with $\eta$ and $I$ in
the sense that
\begin{equation}\label{com}
g(IX,IY)=-g(X,Y),\qquad d\eta(X,Y)=2g(IX,Y), \qquad X,Y\in \mathbb H,
\end{equation}
we have para contact hermitian manifold. The signature of $g$
restricted to $\mathbb H$ is necessarily neutral of type $(n,n)$.

The para-contact Reeb vector field $\xi$ (of length -1) is the
dual vector field to $\eta$ via the metric $g, \quad
g(X,\xi)=\eta(X), \quad \eta(\xi)=-1$ satisfying $d\eta(\xi,.)=0$.
The metric $g$ extends to the metric in the whole manifold by
requiring $g(\xi,\xi)=-1$. In addition, the 1-form $\eta$ is  a
contact form and the fundamental 2-form is defined by
\begin{equation}\label{n1}
2\omega(X,Y)=2g(IX,Y)=d\eta(X,Y).
\end{equation}
The paracomplex structure $I$ on $\mathbb H$ is formally
integrable \cite{Zam} if the $\pm$ eigen-distrubutions $\mathbb
H^{\pm}$ of $I$ in $\mathbb H$ are formally integrable in the
sense that $[\mathbb H^{\pm},\mathbb H^{\pm}]\in \mathbb H^{\pm}
$. Using the Nijenhuis tensor
$N(X,Y)=[IX,IY]+[X,Y]-I[IX,Y]-I[X,IY]$, the formal integrability
of $I$  is equivalent to
\begin{equation}\label{new1}
N(X,Y)=0 \quad {\text and} \quad [IX,Y]+[X,IY]\in \mathbb H.
\end{equation}
A para-contact manifold is called
para-sasakian if $N(X,Y)=d\eta(X,Y)\xi$.

\subsection{The canonical connection}
The canonical  para-contact connection $\nabla$ on a para-contact
hermitian manifold defined in \cite{Zam} is similar to the Webster
connection in the pseudohermitian case. We summarize the
properties of $\nabla$ on an integrable para-contact hermitian
manifold from \cite{Zam}).
\begin{thrm}\cite{Zam} On an integrable para-contact hermitian manifold $(M,\eta,I,g)$ there
exists a unique linear connection preserving  the integrable
para-contact hermitian structure, i.e.
\begin{equation}\label{can1}
\nabla\xi=\nabla I=\nabla\eta=\nabla g=0
\end{equation}
with torsion tensor $T(A,B)=\nabla_AB-\nabla_BA-[A,B]$ given by
\begin{gather}\label{torha}
T(X,Y)=-d\eta(X,Y)\xi=-2\omega(X,Y)\xi, \quad T(\xi,X)\in \mathbb
H,\\\label{cant2}
g(T(\xi,X),Y)=g(T(\xi,Y),X)=g(T(\xi,IX),IY)=\frac12\mathcal
L_{\xi}g(X,Y).
\end{gather}
\end{thrm}
It is shown in \cite{Zam} that the endomorphism $T(\xi,.)$ is the
obstruction an integrable para-contact hermitian manifold to be
parasasakian. We denote the symmetric endomorphism $T_\xi:\mathbb
H\longrightarrow \mathbb H$ by $\tau$ and call it \emph{the
torsion of the integrable para-contact hermitian manifold.} It
follows that the torsion $\tau$ is completely trace-free
\cite{Zam}, i.e.
\begin{equation}\label{tortrace}
\tau(e_a,e_a)=\tau(e_a,Ie_a)=0.
\end{equation}

\section{The Bianchi identities}\label{bib}
Let $R=[\nabla,\nabla]-\nabla_{[,]}$ be the curvature of the
canonical connection $\nabla$. We shall also denote with $R$ the
corresponding (0,4) tensor defined with the help of the metric
$g$. The Ricci tensor $r$, the Ricci 2-form $\rho$ and the
pc-scalar curvature $\text{Scal}$ of $\nabla$ are defined,
respectively, by $$ r(A,B)=R(e_a,A,B,e_a), \qquad
\rho(A,B)=\frac12R(A,B,e_a,Ie_a), \qquad \text{Scal}=r(e_a,e_a),
\quad A,B\in \Gamma(TM).$$
\begin{prop}\label{tech}
Let $(M,\eta,I,g)$ be an integrable para-contact hermitian
manifold. Then:
\begin{itemize}
\item[i)] The curvature of the canonical connection has the properties:
\begin{gather}\label{curi}
R(X,Y,IZ,IV)=-R(X,Y,Z,V), \quad R(X,Y,Z,V)=-R(X,Y,V,Z), \quad
R(X,Y,Z,\xi)=0.
\\\label{currrr}
R(X,Y,Z,V)+R(IX,IY,Z,V)\hspace{8cm}\\\nonumber=
2\Big[g(X,Z)\tau(Y,IV)+g(Y,V)\tau(X,IZ)-g(Y,Z)\tau(X,IV)-g(X,V)\tau(Y,IZ)\Big]\\\nonumber
+2\Big[\omega(X,Z)\tau(Y,V)+\omega(Y,V)\tau(X,Z)-\omega(Y,Z)\tau(X,V)-\omega(X,V)\tau(Y,Z)\Big];\\\label{currr}
R(\xi,X,Y,Z)=(\nabla_Y\tau)(Z,X)-(\nabla_Z\tau)(Y,X).\hspace{6cm}
\end{gather}
\item[ii)] The horizontal Ricci tensor is symmetric,
$r(X,Y)=r(Y,X)$ and has the property
\begin{equation}\label{torric}
r(X,Y)+r(IX,IY)=4(1-n)\tau(X,IY).
\end{equation}
\item[iii)] The horizontal Ricci 2-from satisfies the relations 
\begin{equation}\label{rho}
2\rho(X,IY)=r(X,Y)-r(IX,IY)=R(e_a,Ie_a,X,IY).
\end{equation}
\item[iv)] The following differential identity holds
\begin{equation}\label{div}
2(\nabla_{e_a}r)(e_a,X)=\, d\text{Scal}(X).
\end{equation}
\end{itemize}
\end{prop}
\begin{proof}
Equation \eqref{can1} implies immediately \eqref{curi}. The first
Bianchi identity
\begin{equation}\label{bian1}
\sum_{(A,B,C)}\Bigl\{R(A,B)C-(\nabla_AT)(B,C)-T(T(A,B),C)\Bigr\}=0,
\quad A,B,C\in \Gamma(TM).
\end{equation}
together with \eqref{torha} and \eqref{cant2} yield
\begin{gather}\label{biansim}
R(X,Y,Z,V)-R(Z,V,X,Y)\hspace{9cm}\\\notag\hspace{2cm}=2\omega(X,Z)\tau(Y,V)+2\omega(Y,V)\tau(X,Z)-
2\omega(Y,Z)\tau(X,V)-2\omega(X,V)\tau(Y,Z).
\\\label{biansimv}
R(\xi,X,Y,Z)-R(Y,Z,\xi,X)=(\nabla_Y\tau)(Z,X)-(\nabla_Z\tau)(Y,X).\hspace{4cm}
\end{gather}
Combining \eqref{curi} with \eqref{biansim} and \eqref{biansimv}
we obtain \eqref{currrr} and \eqref{currr}.

When we take the trace of \eqref{biansim}, use \eqref{cant2} and
\eqref{tortrace} we find
$$r(Y,Z)-r(Z,Y)=2\omega(e_a,Z)\tau(Y,e_a)+2\omega(Y,e_a)\tau(e_a,Z)=-2\tau(IZ,Y)+2\tau(Y,IZ)=0.$$
Furthermore, \eqref{curi} and \eqref{currrr} imply
\begin{gather*}r(Y,Z)+r(IY,IZ)=R(e_a,Y,Z,e_a)+R(Ie_a,Y,Z,Ie_a)+4(1-n)\tau(Y,IZ)=4(1-n)\tau(Y,IZ).
\end{gather*}
The first Bianchi identity \eqref{bian1} together with
\eqref{can1} and \eqref{cant2} yields
\begin{equation*}
2\rho(X,IY)=r(X,Y)-r(IY,IX)+
2\tau(X,IY)-2\tau(IY,X)=r(X,Y)-r(IX,IY).
\end{equation*}
The second Bianchi identity reads
\begin{equation}\label{secb}
\sum_{(A,B,C)}\Big\{(\nabla_AR)(B,C,D,E)+R(T(A,B),C,D,E)\Big\}=0,
\qquad A,B,C,D\in\Gamma(TM).
\end{equation}
A suitable trace of  \eqref{secb} leads to
\begin{multline}\label{bian2}
(\nabla_{e_a}R)(X,Y,Z,e_a)-(\nabla_Xr)(Y,Z)+(\nabla_Yr)(X,Z)\\
+2R(\xi,Y,Z,IX)-2R(\xi,X,Z,IY)+2\omega(X,Y)r(\xi,Z)=0.
\end{multline}
The trace  of \eqref{bian2} gives $
2(\nabla_{e_a}r)(X,e_a)-d\text{Scal}(X)+4r(\xi,IX)-4\rho(\xi,X)=0$
while equation \eqref{currr} implies $
r(\xi,IX)=(\nabla_{e_a}\tau)(IX,e_a)=\rho(\xi,X).$ Now, the
identity \eqref{div} follows from  the last two equalities.
\end{proof}

\section{Basic Examples}

Let $\{x_1,y_1,\dots,x_{n+1},y_{n+1}\}$ be the standard coordinate
system in $\mathbb R^{2n+2}$. The standard parahermitian structure
$(\mathbb I,g)$ is defined by
$$\mathbb I\frac{\partial}{\partial x_j}=\frac{\partial}{\partial y_j}, \quad
\mathbb I\frac{\partial}{\partial y_j}=\frac{\partial}{\partial
x_j},\quad g\Big(\frac{\partial}{\partial
x_j},\frac{\partial}{\partial x_k}\Big)=
-g\Big(\frac{\partial}{\partial y_j},\frac{\partial}{\partial
y_k}\Big)=\delta_{jk}, \quad g\Big(\frac{\partial}{\partial
x_j},\frac{\partial}{\partial y_k}\Big)=0,$$ where $\quad j,k
=1,\dots,n$. Recall that a smooth map $f=(u_1,v_1,\dots,u_n,v_n):
\mathbb R^{2n+2}\longrightarrow\mathbb R^{2n+2}$ preserves the
paracomplex structure $\mathbb I$ iff it is \emph{paraholomorphic},
i.e., satisfies the (para) Cauchy-Riemann equations, see
e.g.\cite{Lib}, $df\circ \mathbb I = \mathbb I\circ df$, or,
\begin{equation}\label{PCR}
\frac{\partial u_k}{\partial x_j}=\frac{\partial v_k}{\partial y_j},
\qquad \qquad \frac{\partial u_j}{\partial y_j}=\frac{\partial
v_k}{\partial x_j}.
\end{equation}
Let $(\mathbb R^{2n+2},\mathbb I,g)$ be the standard flat
parahermitian structure on $\mathbb R^{2n+2}$ and $M^{2n+1}$ be a
hypersurface with unit normal $N$ such that $T_pM^{2n+1}\oplus
N=R^{2n+2}, p\in M^{2n+1}$. Consider the vector field $\xi:=\mathbb
IN$, the dual 1-form $\eta(\xi)=-1$ and denote $\mathbb
H=\xi^{\perp}=Ker\,\eta$.   A para CR-structure on $M$ is defined by
$(\mathbb H, I=\mathbb I_{|\mathbb H})$. Moreover
$$d\eta(X,Y)=-\eta([X,Y])=-d\eta(IX,IY)$$ in view of the
integrability condition \eqref{new1}. If in addition
$d\eta_{|\mathbb H}$ is non-degenerate then it necessarily has
signature $(n,n)$ and $(M,\eta)$ is an integrable para-contact
hermitian manifold.
\begin{prop}
Any non-degenerate hypersurface in $(\mathbb R^{2n+2},\mathbb I,g)$
admits an integrable para-contact hermitian structure.
\end{prop}
 Since the horizontal space $\mathbb
H$ is invariant under the standard paracomplex structure of $\mathbb
R^{2n+2}$, a restriction of a paraholomorphic map
$f:R^{2n+2}\longrightarrow R^{2n+2}$ on $(M^{2n+1},\eta)$ induces a para
conformal transformation of the embedded paracontact hermitian structure
$\bar\eta=\mu\eta$ on the hypersurface $M^{2n+1}$.

\subsection{Hyperbolic Heisenberg group}

The hyperbolic Heisenberg group is the  example of an integrable
para-contact hermitian structure with flat canonical connection. The
difference between this group and the standard Heisenberg group is in the
metric, while the groups are identical. {As a group
$\boldsymbol{G\,(\mathbb{P})}\ =\mathbb{R}^{2n}\times \mathbb{R}$ with the
group law given by
\[
( p'',t'')\ =\ (p',t')\circ(p, t)\ =\ \bigl (p'\ +\ p, t' +\ t\ -
\sum_{k=1}^n(u'_kv_k-v'_ku_k)\bigr).
\]
\noindent where $p',p\in\mathbb R^{2n}$, $t', t\in \mathbb R$,
$p=(u_{1},v_{1},\dots,u_{n},v_{n})$ and
$p'=(u'_{1},v'_{1},\dots,u'_{n},v'_{n})$. A basis of left-invariant vector
fields is given by $U_{k}=${\ }$\frac{\partial }{\partial u_{k}}-2v_{k}
\frac{\partial }{\partial t},$ $V_{k}=\frac{\partial }{\partial
v_{k}}+2u_{k} \frac{\partial }{\partial t},\xi =2\frac{\partial }{\partial
t}.$ Define $ \tilde{\Theta}\ =\ -\frac{1}{2}dt\
{-}\sum_{k=1}^{n}(u_{k}dv_{k}-v_{k}du_{k})$ with corresponding horizontal
bundle $\mathbb H$ given by the span of the left-invariant horizontal
vector fields $\{U_{1},...U_{n},V_{1}...V_{n},\}$. We consider an
endomorphism on $\mathbb H$ by defining $IU_{k}${\ }$=V_{k},$ $IV_{k}${\
}$=U_{k} $, hence $I^2=\text{Id}$ on $H$ and $I$ is a paracomplex
structure on $\mathbb H$. The form $\tilde{\Theta}$ and the para-complex
structure $I$ (on $H$) define a para-contact manifold, which will be
called the hyperbolic Heisenberg group. Note that by definition
$\{U_{1},...U_{n},V_{1}...V_{n},\xi \}$ is an orthonormal basis of the
tangent space, $g(U_j,U_j)=-g(V_j,V_j)=1$, $j=1,...,n$.

\begin{thrm}\label{vanh}
Let $(M,\eta,I,g)$ be an integrable para-contact hermitian
manifold of dimension $2n+1$.
\begin{enumerate} \item[i).] If $n>1$ then
$(M,\eta,I,g)$ is locally isomorphic to  the hyperbolic Heisenberg
group exactly when the canonical connection has vanishing
horizontal curvature, $R(X,Y,Z,V)=0$,
\item[ii).] If $n=1$ then $(M,\eta,I,g)$ is locally isomorphic to  the 3-dimensional
hyperbolic Heisenberg group exactly when the canonical connection
has vanishing horizontal curvature and zero torsion.
\end{enumerate}
\end{thrm}
\begin{proof}
It is easy to see that the canonical connection on the hyperbolic
Heisenberg group is the left-invariant connection on the group
which is flat and with zero torsion endomorphism. For the
converse, we first show that if $n>1$ and the horizontal curvature
vanishes then the canonical connection is flat and with zero
torsion endomorphism, $R=\tau=0$. Indeed, \eqref{torric} yields
$\tau=0$ and \eqref{currr} shows $R(\xi,X,Y,Z)=0$.

Let $\{e_1,\dots,e_{n},Ie_1,\dots,Ie_{n},\xi\}$ be a local basis
parallel with respect to $\nabla$. Then \eqref{torha} and
\eqref{cant2} show  that $M$ has the structure of the Lie algebra
of the hyperbolic Heisenberg group, which proves the claim.
\end{proof}

\subsection{Hyperboloid of neutral signature.}\label{standard hyperboloid}

Let $\{x_{0},y_{0},\dots ,x_{n},y_{n}\}$ be the standard coordinate
system in $(\mathbb R^{2n+2},\mathbb I,g)$. Consider the
hypersurface
\begin{equation*}
HS^{2n+1}= \{ (x_{0},y_{0},\dots ,x_{n},y_{n})\subset \mathbb
R^{2n+2}\, |\ x_{0}^{2}+\dots +x_{n}^{2}-y_{0}^{2}\dots -y_{n}^{2}=1
\}.
\end{equation*}
$HS^{2n+1}$ carries a natural para-CR structure inherited from its
embedding in $(\mathbb R^{2n+2},\mathbb I,g)$. The horizontal bundle
$\mathbb{H}$ is the maximal subspace of the tangent space of
$HS^{2n+1}$ which is invariant under the (restriction of the) action
of $\mathbb I$. We take
\begin{equation*}
\tilde\eta =-\sum_{j=0}^{n}\left( x_{j}dy_{j}-y_{j}dx_{j}\right)
\end{equation*}
noting that here $N=\sum_{j=0}^{n}\left( x_{j}\frac{\partial }{\partial x_{j}
}+y_{j}\frac{\partial }{\partial y_{j}}\right) $ and $\xi
=\sum_{j=0}^{n}\left( x_{j}\frac{\partial }{\partial y_{j}}+y_{j}\frac{
\partial }{\partial x_{j}}\right) .$  We will also  consider  $HS^{2n+1}$ as the boundary of
the "ball"
 $B= \{ (x_{0},y_{0},\dots ,x_{n},y_{n})\subset
\mathbb R^{2n+2}\, : \ x_{0}^{2}+\dots +x_{n}^{2}-y_{0}^{2}\dots
-y_{n}^{2} < 1 \}.$

\subsection{The Cayley transform}
Let $\Sigma_0=\{(x_{0},y_{0},\dots ,x_{n},y_{n})\in HS^{2n+1} :
(1+x_0)^2=y_0^2 \}$. The
 Cayley transform (centered at $\Sigma_0$), is defined as follows
\begin{equation}\label{e:cayley}
\begin{aligned}
& \hskip1truein \mathcal{C}: HS^{2n+1}\setminus \Sigma_0\rightarrow
\boldsymbol{
G\,(\mathbb{P})}\\
&  t=\frac {2y_0}{(1+x_0)^2-y_0^2} ,\qquad u_k =\frac
{x_k(1+x_0)-y_ky_0}{(1+x_0)^2-y_0^2}, \qquad v_k =\frac
{y_k(1+x_0)-x_ky_0}{(1+x_0)^2-y_0^2}.
\end{aligned}
\end{equation}
 A small
calculation shows
\[
\mathcal{C}^*\tilde{\Theta} = \frac {1}{(1+x_0)^2-y_0^2} \tilde\eta.
\]
Furthermore, the para-complex structure is preserved. In order to
see the last claim, we can consider $\GP$ as the boundary of the
domain $D=\{ (u_{0},v_{0},\dots ,u_{n},v_{n})\subset \mathbb
R^{2n+2}\, : \ u_{1}^{2}+\dots +u_{n}^{2}-v_{1}^{2}\dots
-v_{n}^{2}< v_0 \}$ by identifying the point $(p,t)\in\GP$ with
the point $(t,\Sigma_{k=1}^n (u_k^2-v_k^2),u_{1},v_{1}\dots
,u_{n},v_{n})\in\partial D$ and define the diffeomorphism $
\mathcal{C}:B\setminus \Sigma_0\rightarrow D\setminus\Xi_0$,
$\Xi_0= \{ (1+u_0)^2-v_0^2=0 \}$,
\begin{equation*}
\begin{aligned}
&  u_0=\frac {2y_0}{(1+x_0)^2-y_0^2} ,\qquad  v_0=\frac {1-x_0^2+y_0^2}{(1+x_0)^2-y_0^2} \\
& u_k =\frac {x_k(1+x_0)-y_ky_0}{(1+x_0)^2-y_0^2}, \qquad v_k =\frac
{y_k(1+x_0)-x_ky_0}{(1+x_0)^2-y_0^2}.
\end{aligned}
\end{equation*}
A calculation shows that the above map is para-holomorphic, i.e.,
the coordinate functions satisfy the (para) Cauchy-Riemann equations
\eqref{PCR}.  Thus, $\mathcal{C}$ preserves the para-CR structure
when considered as a map between the boundaries of $B$ and $D$.

Using hyperbolic rotations, which preserve the para-contact
structure, and Cayley maps similar to the above  we see that the
hyperboloid is locally para-contact conformal to the hyperbolic
Heisenberg group.

{ Finally, it is worth noting that according to Theorem \ref{main2} the
hyperboloid $HS^{4n+1}$is para-contact conformally flat, while regarded as
a CR submanifold of $\mathbb{C}^{2n+1}$ it is  a degenerate CR manifold}.

\section{Paracontact conformal curvature}

In this section we define para-contact conformal invariant and
prove Theorem~\ref{main1}.

\subsection{Paracontact conformal transformations}\label{s:conf transf}

A conformal para-contact transformation (\emph{pc transformartion}
) between two para-contact manifold is a diffeomorphism $\Phi$
which preserves the para-contact structure i.e.
$
\Phi^*\eta=\mu\eta,
$ for a nowhere vanishing smooth function $\mu$.

Let $u$ be a smooth nowhere vanishing function on a para-contact
manifold $(M, \eta)$. Let $\bar\eta=\frac{1}{2}e^{-2u}\eta$ be a
conformal deformation of  $\eta$. We will denote the objects related
to $\bar\eta$ by over-lining the same object corresponding to
$\eta$. Thus,
$
 d\bar\eta=-e^{-2u}du\wedge\eta\ +\ \frac{1}{2}e^{-2u}d\eta,\qquad
\bar g=\frac{1}{2}e^{-2u}g. $ The new  para-contact Reeb vector
field  $\bar\xi$ is \cite{Zam}
\begin{equation}\label{New19}
\bar\xi\ =\ 2e^{2u}\,\xi\ + 2e^{2u}\ I\nabla u,
\end{equation}
where
$\nabla u$ is the horizontal gradient defined by $g(\nabla u,X)=du(X),\quad X\in
H$.
The horizontal sub-Laplacian and the norm of the horizontal gradient  are defined
respectively by
$\triangle u\ =\ tr^g_H(\nabla du)\  = \ \nabla du(e_a,e_a)=\sum_{s=1}^n(\nabla du(e_s,e_s)-\nabla du(Ie_s,Ie_s)),\qquad
 |\nabla u|^2\ =\ du(e_a)^2=\sum_{s=1}^n(du(e_s)^2-du(Ie_s)^2).
$ The canonical para-contact connections $\nabla$ and $\bar\nabla$
are related by a (1,2) tensor S,
\begin{equation}\label{qcw2}
\bar\nabla_AB=\nabla_AB+S_AB, \qquad A,B\in\Gamma(TM).
\end{equation}
Suppose the para contact structure is integrable. The conditions
\eqref{torha} and $\bar\nabla\bar g=0$ determine $g(S(X,Y),Z)$ for
$X,Y,Z\in H$ due to the equality
\begin{multline}\label{Ivan4}
g(S(X,Y),Z)=-du(X)g(Y,Z)-du(IX)\omega(Y,Z)\\-
du(Y)g(Z,X)+du(IY)\omega(Z,X)+du(Z)g(X,Y)+du(IZ)\omega(X,Y).
\end{multline}
We obtain after some calculations using \eqref{New19}  that
\begin{equation}\label{New20}
\bar\tau(X,Y)-2e^{2u}\tau(X,Y)-g(S(\bar\xi,X),Y)=-2e^{2u}\nabla
du(X,IY)-4e^{2u}du(X)du(IY).
\end{equation}
From \eqref{New20} and \eqref{cant2} we find
\begin{multline}\label{New220}
g(S(\bar\xi,X),Y)-g(S(\bar\xi,IX)IY)\\
=2e^{2u}\Bigl[\nabla du(X,IY)-\nabla
du(IX,Y)+2du(X)du(IY)-2du(IX)du(Y)\Bigr].
\end{multline}
The condition $\bar\nabla I=\nabla I=0$ yield
$g(S(\bar\xi,X),Y)=-g(S(\bar\xi,IX)IY)$. Substitute the latter
into \eqref{New20} and \eqref{New220}, use \eqref{New19} and
\eqref{Ivan4} to get
\begin{gather}\label{qcw3}
 g(S(\xi,X),Y)=\frac{1}{2}\Big[\nabla du(X,IY)-\nabla du(IX,Y)\Big]\hspace{6.5cm}\\\nonumber \hspace{7cm}
 -du(X)du(IY)+du(IX)du(Y)
 +|\nabla u|^2\omega(X,Y),\\\label{contor}
 \bar\tau(X,Y)=e^{2u}\Big[2\tau(X,Y)-\nabla du(X,IY)-\nabla du(IX,Y)-2du(X)du(IY)-2du(IX)du(Y)\Big].
\end{gather}
\noindent In addition, the pc-scalar curvature changes according to the formula
\cite{Zam}
\begin{equation}\label{e:conf change scalar curv}
\overline {\text{Scal}}\ =\ 2e^{2u}\,\text{Scal}\ -\
8n(n+1)e^{2u}|\nabla u|^2\ +8(n+1)e^{2u}\triangle u.
\end{equation}
The identity $d^2=0$ yields
$\nabla du(X,Y)-\nabla du(Y,X)=-du(T(X,Y)).$
Applying \eqref{torha}, we can write
\begin{equation}\label{symdh}
\nabla du(X,Y)=[\nabla du]_{[sym]}(X,Y)+du(\xi)\omega(X,Y),
\end{equation}
where $[.]_{[sym]}$ denotes the symmetric part of the corresponding (0,2)-tensor.

\subsection{Paracontact conformal curvature tensor}
Let $(M,\eta,I,g)$ be a (2n+1)-dimensional integrable para-contact
hermitian manifold. Let us consider the symmetric (0,2) tensor $L$
defined on $H$ by the equality
\begin{equation}\label{lll}
L(X,Y)=
\frac1{2(n+2)}\rho(X,IY)-\tau(IX,Y)-\frac{\text{Scal}}{8(n+1)(n+2)}g(X,Y),
\qquad X,Y\in H.
\end{equation}
We define the (0,4) tensor $PW$ on $H$ by
\begin{multline}\label{qcwdef}
g(PW(X,Y)Z,V)= g(R(X,Y)Z,V)\\
+g(X,Z)L(Y,V)+g(Y,V)L(X,Z)-g(Y,Z)L(X,V)-g(X,V)L(Y,Z)\\
+\omega(X,Z)L(Y,IV)+\omega(Y,V)L(X,IZ)-\omega(Y,Z)L(X,IV)-\omega(X,V)L(Y,IZ)\\
+\omega(X,Y)\Bigl[L(Z,IV)-L(IZ,V) \Bigr]
+\omega(Z,V)\Bigl[L(X,IY)-L(IX,Y)\Bigr].
\end{multline}

\begin{prop}\label{trfree}
The tensor $PW$ is completely trace-free, i.e.
$$r(PW)=\rho(PW)=0.
$$
\end{prop}
\begin{proof}
The claim follows after taking the corresponding traces in
\eqref{qcwdef} keeping in mind \eqref{lll}.
\end{proof}
If we compare \eqref{qcwdef} and \eqref{currr} we obtain the
following
\begin{prop}\label{main0}
For $n>1$ the  tensor $PW$ has the properties
$$PW(X,Y,Z,V)+PW(IX,IY,Z,V)=0,$$
\begin{multline}\label{pinv}
PW(X,Y,Z,V)-PW(IX,IY,Z,V)=R(X,Y,Z,V)-R(IX,IY,Z,V)\\
+\frac{\text{Scal}}{2(n+1)(n+2)}\Bigl[\omega(X,Z)\omega(Y,V)
-\omega(Y,Z)\omega(X,V)+2\omega(X,Y)\omega_s(Z,V)
 \Bigr]\\ -\frac{\text{Scal}}{2(n+1)(n+2)}\Bigl[g(X,Z)g(Y,V)-g(Y,Z)g(X,V)\Bigr]
 +\frac2{n+2}\Big[\omega(X,Y)\rho(Z,V)+\omega(Z,V)\rho(X,Y)\Big]\\
 +\frac1{n+2}\Big[g(X,Z)\rho(Y,IV)-g(Y,Z)\rho(X,IV)+g(Y,V)\rho(X,IZ)-g(X,V)\rho(Y,IZ)\Big]\\
 +\frac1{n+2}\Big[\omega(X,Z)\rho(Y,V)-\omega(Y,Z)\rho(X,V)+\omega(Y,V)\rho(X,Z)-\omega(X,V)\rho(Y,Z)\Big]\\.
\end{multline}
For $n=1$ the tensor $PW$ vanishes identically.
\end{prop}
\begin{dfn}\label{d:Wpc}
We denote  the tensor $PW(X,Y,Z,V)-PW(IX,IY,Z,V)$ by $2W^{pc}$ and
call it \emph{the para-contact conformal curvature}.

If $n=1$ we define on $\mathbb H$ the following symmetric (0,2)
tensor $F$ by the equality
\begin{gather}\label{treeF}
F(X,Y)=(\nabla d(Scal))(X,IY)+(\nabla d(Scal))(Y,IX)
+16(\nabla^2_{Xe_a}\tau)(Y,e_a)+16(\nabla^2_{Ye_a}\tau)(X,e_a)
\\\nonumber \hspace{1.5cm}
-48(\nabla^2_{e_aIe_a}\tau)(X,IY)+36Scal\cdot\tau(X,Y)+3g(X,Y)(\nabla
d(Scal))(e_aIe_a).
\end{gather}
\end{dfn}

\subsection{Proof of Theorem~\ref{main1}} First we show
\begin{thrm}\label{qcinv}
The para-contact conformal curvature $W^{pc}$ of an integrable
para-contact hermitian manifold is invariant under conformal
para-contact transformations, i.e., if
$$2\bar\eta=e^{-2u}\eta\quad {\text for}\quad {\text any}\quad
{\text smooth}\quad {\text function}\quad u \quad {\text then}
\qquad 2e^{2u}W^{pc}_{\bar\eta}=W^{pc}_{\eta}.$$
\end{thrm}
\begin{proof}
After a  straightforward computation using \eqref{qcw2},
\eqref{Ivan4} and \eqref{qcw3} we obtain the formula
\begin{multline}\label{qcw4}
2e^{2u}g(\bar R(X,Y)Z,V)-g(R(X,Y)Z,V)=-g(Z,V)\Bigl[M(X,Y)-M(Y,X)\Bigr]\\
-g(X,Z)M(Y,V)-g(Y,V)M(X,Z)+g(Y,Z)M(X,V)+g(X,V)M(Y,Z)\\
-\omega(X,Z)M(Y,IV)-\omega(Y,V)M(X,IZ)+\omega(Y,Z)M(X,IV)+\omega(X,V)M(Y,IZ)\\
-\omega(X,Y)\Bigl[M(Z,IV)-M(IZ,V)\Bigr]
-\omega(Z,V)\Bigl[M(X,IY)-M(Y,IX)\Bigr].
\end{multline}
where the (0,2) tensor $M$ is given by
\begin{equation}\label{qcw5}
M(X,Y)=\nabla du(X,Y)+du(X)du(Y)+ du(IX)du(IY)-\frac12g(X,Y)|\nabla
u|^2.
\end{equation}
Let $tr M=M(e_a,e_a)$ be the trace of the tensor $M$. Using
\eqref{qcw5} and \eqref{symdh} we obtain
\begin{gather}\label{qcw6}
tr M=\triangle u-n|\nabla u|^2,\quad
M(X,Y)+M(IX,IY)=M(Y,X)+M(IY,IX), \\\label{qcw6a}
M(e_a,Ie_a)=-2ndu(\xi), \quad
M(e_a,Ie_a)\omega(X,Y)=-n\Big[M(X,Y)-M(Y,X)\Big].
\end{gather}
Taking the trace in \eqref{qcw4} and using \eqref{qcw5},
\eqref{qcw6}, and \eqref{qcw6a} we come to
\begin{gather}\label{qcwric}
\overline
r(X,Y)-r(X,Y)\hspace{11cm}\\\nonumber=(n+1)M(X,Y)+nM(Y,X)-M(IX,IY)-2M(IY,IX)+tr
M\,g(X,Y);\\\nonumber
e^{-2u}\overline{Scal}-2\,\text{Scal}=8(n+1)tr M.\hspace{7cm}
\end{gather}
Proposition~\ref{tech} together with \eqref{qcwric} and \eqref{lll} imply
\begin{equation}\label{mm}
M_{[sym]}(X,Y)= \overline L(X,Y)-L(X,Y).
\end{equation}
Now, from \eqref{qcw5} and \eqref{symdh} we obtain
\begin{equation}\label{mm1}
M(X,Y)=M_{[sym]}(X,Y)+du(\xi)\omega(X,Y).
\end{equation}
Substituting \eqref{mm} into \eqref{mm1}, then inserting the obtained
equality in \eqref{qcw4} and finally using \eqref{qcw6}  allows us to
complete the proof of Theorem~\ref{qcinv}.
\end{proof}

At this point a combination of Theorem~\ref{qcinv} and
Proposition~\ref{main0} ends the proof of Theorem~\ref{main1}.

\section{Converse problem. Proof of Theorem~\ref{main2}}

Suppose  $W^{pc}=0$, hence   $PW=0$ by Proposition~\ref{main0}. We
shall show that in this case there  exists (locally) a smooth
conformal factor $u$ which changes by a pc conformal
transformation the integrable para-contact hermitian structure to
a torsion-free flat one.

Consider  the following system of differential equations with
respect to the unknown function $u$
\begin{gather}\label{sist1}
\nabla du(X,Y)=-L(X,Y)-du(X)du(Y)-du(IX)du(IY)+ \frac12g(X,Y)|\nabla
u|^2 + du(\xi)\omega(X,Y)\\\label{add1} \nabla du(X,\xi)=-\mathbb
B(X,\xi)-L(X,I\nabla u)+\frac12du(IX)|\nabla u|^2
-du(X)du(\xi_i)\hspace{3.5cm}\\\label{add2} \nabla
du(\xi,\xi)=-\mathbb B(\xi,\xi) -\mathbb B(I\nabla
u,\xi)-\frac14|\nabla u|^4-(du(\xi))^2,\hspace{5.5cm}
\end{gather}
where  $\mathbb B(X,\xi)$ and $\mathbb B(\xi,\xi)$ do not depend on the function $u$
and are determined in \eqref{bes} and
\eqref{bst}.

In order to prove Theorem~\ref{main2} it is sufficient to show the
existence of a local smooth solution to \eqref{sist1}. Indeed,
suppose $u$ is a local smooth solution to \eqref{sist1}. Then the
canonical connection of the para-contact hermitian structure
$2\bar\eta=e^{-2u}\eta$ has in view of \eqref{contor} zero torsion.
Furthermore,  the curvature restricted to $H$ vanishes when
$W^{pc}=0$ taking into account Proposition~\ref{main0} and the proof
of Theorem~\ref{qcinv}. Therefore, we can apply Theorem~\ref{vanh}
to conclude the result.

The rest of this section is devoted to showing the existence of a
smooth solution to the system \eqref{sist1}-\eqref{add2}.

The integrability conditions for this overdetermined system are the Ricci identities,
\begin{equation}\label{integr}
\nabla du(A,B,C)-\nabla du(B,A,C)=-R(A,B,C,\nabla u)-\nabla
du((T(A,B),C), \quad A,B,C \in \Gamma(TM).
\end{equation}

We  consider as separate cases the four possibilities for $A$, $B$
and $C$.

{\bf Case 1: $\Bigl[A,B,C \in H\Bigr]$}.  Invoking \eqref{torha} we see
that equation \eqref{integr} on $H$ has the  form
\begin{equation}\label{inteh}
\nabla du(Z,X,Y)-\nabla du(X,Z,Y)+R(Z,X,Y,\nabla
u)-2\omega(Z,X)\nabla du(\xi,Y)=0,
\end{equation}

Take a covariant derivative of \eqref{sist1} along $Z\in H$, substitute in
the obtained equality \eqref{sist1} and \eqref{lll}, anticommute the
covariant derivatives, let $W^{pc}=0$ in \eqref{qcwdef}, and finally use
\eqref{add1} and \eqref{sist1} to see that the integrability condition
\eqref{inteh} is
\begin{equation}\label{inte}
(\nabla_ZL)(X,Y)-(\nabla_XL)(Z,Y)=\omega(Z,Y)\mathbb B(X,\xi)-
\omega(X,Y)\mathbb B(Z,\xi)+2\omega(Z,X)\mathbb B(Y,\xi).
\end{equation}
The 1-forms $\mathbb B(X,\xi)$ can be determined by taking traces in
\eqref{inte}. Thus, we have
\begin{equation}\label{bes}
(\nabla_{e_a}L)(Ie_a,IX)=-(2n+1)\mathbb
B(IX,\xi)\quad\text{and}\quad
(\nabla_Xtr\,L)-(\nabla_{e_a}L)(e_a,X)=3\mathbb B(IX,\xi).
\end{equation}
Notice that the consistence of the first and second equalities in
\eqref{bes} is  equivalent to \eqref{div}.
\begin{lemma}\label{integr1}
Suppose $W^{pc}=0$ and the dimension  is bigger than three. Then
\eqref{inte} holds.
\end{lemma}
\begin{proof}
Using \eqref{torha}, the second Bianchi  identity \eqref{secb} gives
\begin{multline}\label{bi2ric}
(\nabla_{e_a}R)(X,Y,Z,e_a)-(\nabla_Xr)(Y,Z)+(\nabla_Yr)(X,Z)\\
+2R(\xi,Y,Z,IX)-2R(\xi,X,Z,IY)+2\omega(X,Y)r(\xi,Z)=0,
\end{multline}
\begin{multline}\label{mon20}
(\nabla_{X}\rho)(Y,Z)+(\nabla_{Y}\rho)(Z,X)+(\nabla_{Z}\rho)(X,Y)-\\
-2\omega(X,Y)\rho(\xi,Z)-2\omega(Y,Z)\rho(\xi,X)-2\omega(Z,X)\rho(\xi,Y)=0,
\end{multline}
\begin{equation}\label{mon27}
(\nabla_{X}\rho)(Y,Z)+(\nabla_{e_a}R)(Ie_a,X,Y,Z)+2(n-1)R(\xi,X,Y,Z)=0.
\end{equation}
From $W^{pc}=0$ and \eqref{lll} we can express $r, \rho$ and $\tau$ in
terms of $L$ and $tr\,L$, obtaing
\begin{align}\label{rl0}
r(X,Y) & =(2n+1)L(X,Y)-3L(IX,IY)+(trL)g(X,Y)\\\label{rl1}
\rho(X,Y)& =(n+2)L(X,IY)-(n+2)L(IX,Y)-(trL)\omega(X,Y)\\\label{tl}
2\tau(IX,Y) & =-L(X,Y)-L(IX,IY).
\end{align}
Inserting \eqref{qcwdef} and \eqref{currr} in \eqref{bi2ric}, and
then  using \eqref{rl0}, \eqref{tl} we come after some standard
calculations to the following identity
\begin{multline}\label{mon19}
-3g(Z,X)\mathbb B(IY,\xi)+3g(Z,Y)\mathbb
B(IX,\xi)-(2n+1)\omega(X,Z)\mathbb B(Y,\xi)\\
+(2n+1)\omega(Y,Z)\mathbb B(X,\xi)-2(2n+1)\omega(X,Y)\mathbb
B(Z,\xi)+2n[(\nabla_XL)(Y,Z)-(\nabla_YL)(X,Z)]\\
+[(\nabla_{IZ}L)(X,IY)-(\nabla_{IZ}L)(IX,Y)]-[(\nabla_{IX}L)(IY,Z)-(\nabla_{IY}L)(IX,Z)]\\
-2[(\nabla_{IX}L)(Y,IZ)-(\nabla_{IY}L)(X,IZ)]-3[(\nabla_XL)(IY,IZ)-(\nabla_YL)(IX,IZ)]=0.
\end{multline}
A substitution  of \eqref{qcwdef} and \eqref{currr} in \eqref{mon20}
together with \eqref{rl1} give
\begin{multline}\label{mon21}
(n+2)[(\nabla_XL)(Y,IZ)-(\nabla_YL)(X,IZ)]-(n+2)[(\nabla_XL)(IY,Z)-(\nabla_YL)(IX,Z)]\\
+(n+2)[(\nabla_ZL)(X,IY)-(\nabla_ZL)(IX,Y)]\\-2(n+2)[\omega(X,Y)\mathbb
B(IZ,\xi)+\omega(Y,Z)\mathbb B(IX,\xi)+\omega(Z,X)\mathbb
B(IY,\xi)]=0.
\end{multline}
Take $ IZ$ instead of  $Z$ in \eqref{mon21}, then set $IX$ and $ IY$,
correspondingly, for $X$ and $Y$ into the obtained result. Taking the sum
of thus achieved equalities we derive
\begin{multline}\label{mon24}
[(\nabla_XL)(Y,Z)-(\nabla_YL)(X,Z)]-[(\nabla_XL)(IY,IZ)-(\nabla_YL)(IX,IZ)]+\\
+[(\nabla_{IX}L)(IY,Z)-(\nabla_{IY}L)(IX,Z)]-[(\nabla_{IX}L)(Y,IZ)-(\nabla_{IY}L)(X,IZ)]+\\
+2g(Y,Z)\mathbb B(IX,\xi)-2g(Z,X)\mathbb
B(IY,\xi)+2\omega(Y,Z)\mathbb B(X,\xi)-2\omega(X,Z)\mathbb
B(Y,\xi)=0.
\end{multline}
Insert \eqref{qcwdef}, \eqref{currr} in \eqref{mon27} using \eqref{rl1},
\eqref{tl}, replace $Y$ and $Z$ respectively with $ IY$ and $ IZ$ into the
obtained equality  and then take the sum of both equations to obtain
\begin{multline}\label{mon30}
(n-1)[(\nabla_{IZ}L)(X,IY)-(\nabla_{IY}L)(X,IZ)]+(n-1)[(\nabla_{IZ}L)(IX,Y)-(\nabla_{IY}L)(IX,Z)]\\
+(n-1)[(\nabla_{Z}L)(X,Y)-(\nabla_{Y}L)(X,Z)]+(n-1)[(\nabla_{Z}L)(IX,IY)-(\nabla_{Y}L)(IX,IZ)]=0.
\end{multline}
Substitute $X$ by $Z$, and $Z$ by $X$ in \eqref{mon30}. The sum of the
obtained equalities and \eqref{mon24} yield
\begin{multline}\label{mon32}
[(\nabla_XL)(Y,Z)-(\nabla_YL)(X,Z)]+[(\nabla_{IX}L)(IY,Z)-(\nabla_{IY}L)(IX,Z)]-\\
 -g(X,Z)\mathbb B(IY,\xi)+g(Y,Z)\mathbb
B(IX,\xi)-\omega(X,Z)\mathbb B(Y,\xi)+\omega(Y,Z)\mathbb
 B(X,\xi)=0.
\end{multline}
The cyclic sum in \eqref{mon32} gives
\begin{multline}\label{mon33}
[(\nabla_{IZ}L)(X,IY)-(\nabla_{IZ}L)(IX,Y)]=[(\nabla_{IX}L)(IY,Z)-(\nabla_{IY}L)(IX,Z)]-\\
-[(\nabla_{IX}L)(Y,IZ)-(\nabla_{IY}L)(X,IZ)]+2\omega(Z,X)\mathbb
B(Y,\xi)+2\omega(Y,Z)\mathbb B(X,\xi)+2\omega(X,Y)\mathbb B(Z,\xi).
\end{multline}
Now,  identity \eqref{inte} follows from \eqref{mon19}, \eqref{mon32} and \eqref{mon33}.
\end{proof}

{\bf Case 2: $\Bigl[A,B \in H,\quad C=\xi\Bigr]$}. In this case, with the
help of \eqref{torha},  \eqref{integr} turns into the equation
\begin{multline}\label{intehxi}
\nabla du(Z,X,\xi)-\nabla du(X,Z,\xi)=-R(Z,X,\xi,\nabla u)-\nabla
du(T(Z,X),\xi)= 2\omega(Z,X)\nabla du(\xi,\xi).
\end{multline}
Take a covariant derivative of \eqref{add1} along $Z\in H$, substitute
\eqref{sist1} and \eqref{add1} in the obtained equality, then anticommute
the covariant derivatives and substitute the result in \eqref{intehxi}
together with the already established \eqref{inte}, \eqref{add2} and
\eqref{lll} to get after some standard calculations that the integrability
condition in this case is
\begin{equation}\label{inte1}
(\nabla_Z\mathbb B)(X,\xi)-(\nabla_X\mathbb
B)(Z,\xi)+L(Z,IL(X))-L(X,IL(Z))= 2\mathbb B(\xi,\xi)\omega(Z,X).
\end{equation}
Here, the function $\mathbb B(\xi,\xi)$ is independent of $u$ and is
uniquely determined  by
\begin{equation}\label{bst}
\mathbb B(\xi,\xi)=-\frac{1}{2n}[(\nabla_{e_a}\mathbb
B)(Ie_a,\xi)+L(e_a,IL(Ie_a))].
\end{equation}
\begin{lemma}\label{integr2}
If $W^{pc}=0$ and the dimension is bigger than three, then
\eqref{inte1} holds.
\end{lemma}
\begin{proof}
Differentiate the already proved \eqref{inte}, take the
corresponding traces and use the symmetry of $L$ to see
\begin{align}\label{pp1}%
& (\nabla^2_{e_a,Ie_a}L)(Y,Z)-(\nabla^2_{e_a,Y}L)(Ie_a,Z)  \\\notag
& \hskip1truein =(\nabla_Z \mathbb
B)(Y,\xi)-\omega(Y,Z)(\nabla_{e_a} \mathbb B)(Ie_a,\xi)+2(\nabla_Y
\mathbb B)(Z,\xi) \\\label{pp3} &
-(\nabla^2_{e_a,Y}L)(Ie_a,Z)+(\nabla^2_{e_a,Z}L)(Ie_a,Y)
\\\notag & \hskip1truein = -(\nabla_Z \mathbb B)(Y,\xi)-2\omega(Y,Z)(\nabla_{e_a}
\mathbb B)(Ie_a,\xi)+(\nabla_Y \mathbb B)(Z,\xi)
\\\label{pp2}
& (\nabla^2_{Y,e_a}L)(Ie_a,Z)=-(2n+1)(\nabla_Y \mathbb B)(Z,\xi).
\end{align}
A combination of \eqref{pp1}, \eqref{pp2} and  \eqref{pp3} yields
\begin{multline}\label{pp4}
\Big[(\nabla^2_{Y,e_a}L)-(\nabla^2_{e_a,Y}L)\Big](Ie_a,Z)-
\Big[(\nabla^2_{Z,e_a}L)-(\nabla^2_{e_a,Z}L)\Big](Ie_a,Y)\\=2n(\nabla_Z
\mathbb B)(Y,\xi)-\omega(Y,Z)(\nabla_{e_a} \mathbb
B)(Ie_a,\xi)-2n(\nabla_Y \mathbb B)(Z,\xi).
\end{multline}
The Ricci identities, \eqref{torha}, \eqref{rho}, \eqref{rl0},
\eqref{rl1} and \eqref{tl} give
\begin{multline}\label{pp5}
\Big[(\nabla^2_{Y,e_a}L)-(\nabla^2_{e_a,Y}L)\Big](Ie_a,Z)=2(\nabla_{\xi}L)(Y,Z)-(trL)[L(Y,IZ)+L(IY,Z)]\\
+(2n+1)L(Y,IL(Z))-3L(IY,L(Z))
-3L(IZ,L(Y))+L(Z,IL(Y))+\omega(Y,Z)L(e_a,IL(e_a)).
\end{multline}
\begin{multline}\label{pp6}
(\nabla^2_{e_a,Ie_a}L)(Y,Z)=
(n+2)[L(IY,L(Z))-L(Y,IL(Z))]\\
-(n+2)L(Z,IL(Y))+(n+2)L(IZ,L(Y))-2n(\nabla_{\xi}L)(Y,Z)+(trL)(L(IY,Z)+L(Y,IZ)).
\end{multline}
The identity \eqref{inte1} follows from \eqref{pp4} and \eqref{pp5}.
\end{proof}

{\bf Case 3: $\Bigl[A=\xi, \quad B, C \in H\Bigr]$}. In this case
\eqref{integr} becomes
\begin{multline}\label{intexih1}
\nabla du(\xi,X,Y)-\nabla du(X,\xi,Y)+R(\xi,X,Y,\nabla u)+\nabla
du(T(\xi,X),Y)=0.\hfill
\end{multline}
Take the covariant derivative of \eqref{sist1} along $\xi$ and a covariant
derivative of \eqref{add1} along a horizontal direction, apply
\eqref{add1}, \eqref{sist1}, \eqref{add2}, use \eqref{currr} and a
suitable traces of \eqref{sist1} and \eqref{tl} to get from
\eqref{intexih1} with the help of \eqref{tl}, \eqref{lll} and the already
proved \eqref{inte} that the integrability condition \eqref{intexih1}
becomes
\begin{multline}\label{intexih11}
(\nabla_X\mathbb
B)(Y,\xi)-(\nabla_{\xi}L)(X,Y)=L(Y,IL(X))+\tau(X,L(Y))+\tau(Y,L(X))+\mathbb
B(\xi,\xi)\omega(X,Y).
\end{multline}
Notice that Case 3 implies Case 2 since  \eqref{inte1} is the
skew-symmetric part of \eqref{intexih11}.
\begin{lemma}\label{integr3}
Suppose $W^{pc}=0$ and dimension is bigger than 3. Then \eqref{intexih11} holds.
\end{lemma}
\begin{proof}
Combine \eqref{pp1}, \eqref{pp2}, \eqref{pp3} and the already proved \eqref{inte1} to obtain
\begin{multline}\label{pp7}
(\nabla^2_{e_a,Ie_a}L)(Y,Z)+\Big[(\nabla^2_{Y,e_a}L)-(\nabla^2_{e_a,Y}L)\Big](Ie_a,Z)=-2(n-1)(\nabla_Y
\mathbb B)(Z,\xi)-\\
-2\omega(Y,Z)\mathbb
B(\xi,\xi)+L(Y,IL(Z))-L(Z,IL(Y))-\omega(Y,Z)(\nabla_{e_a} \mathbb
B)(Ie_a,\xi)
\end{multline}
Now, \eqref{pp5}, \eqref{pp6} and \eqref{pp7} imply \eqref{intexih11}.
\end{proof}

{\bf Case $4$: $\Bigl[A \in H,\quad B=C=\xi\Bigr]$}. In this case
\eqref{integr} has the form
\begin{equation}\label{intehxi32}
\nabla du(X,\xi,\xi)-\nabla du(\xi,X,\xi)=-R(X,\xi,\xi,\nabla
u)+\nabla du(T(\xi,X),\xi)=\tau(X,e_a)\nabla du(e_a,\xi).
\end{equation}
Take the covariant derivative of \eqref{add1} along  $\xi$  and a
covariant derivative of \eqref{add2} along a horizontal direction, then
use \eqref{sist1}, the already proved  \eqref{intexih11}, apply
\eqref{add1} to see that \eqref{intehxi32} is equivalent to
\begin{equation}\label{intehxi312}
(\nabla_{\xi}\mathbb B)(X,\xi)-(\nabla_X\mathbb
B)(\xi,\xi)-2\mathbb B(e_a,\xi)L(X,Ie_a)+\tau(X,e_a)\mathbb
B(e_a,\xi)=0.
\end{equation}
\begin{lemma}\label{integr4}
Suppose $W^{pc}=0$ and dimension is bigger than 3. Then \eqref{intehxi312} holds.
\end{lemma}
\begin{proof}
Differentiate the already proven \eqref{inte1}, \eqref{intexih11},
the first equality in \eqref{bes}, take the corresponding traces,
 use the symmetry of $L$, $\tau$ and \eqref{bes} to get
\begin{gather}\label{pp41}
(\nabla^2_{e_a,Ie_a}\mathbb
B)(Y,\xi)=(n+2)[L(IY,e_b)-L(Y,Ie_b)]\mathbb
B(e_b,\xi)-(trL)\mathbb B(IY,\xi)-2n(\nabla_{\xi}\mathbb B)(Y,\xi)
\\\label{pp42}
(\nabla^2_{e_b,Ie_b}\mathbb
B)(Y,\xi)-(\nabla^2_{e_b,\xi}L)(Ie_b,Y)=-(2n+1)\mathbb
B(e_a,\xi)[L(Y,Ie_a)+\tau(Y,e_a)] +(\nabla_{Y}\mathbb
B)(\xi,\xi)\\\nonumber+[(\nabla_{e_b}L)(Y,Ie_a)+(\nabla_{e_b}\tau)(Y,e_a)]L(Ie_b,e_a)+
(\nabla_{e_b}\tau)(Ie_b,e_a)L(Y,e_a)
+\tau(Ie_b,e_a)(\nabla_{e_b}L)(Y,e_a)
\\\label{pp43}
(\nabla^2_{\xi,e_b}L)(Ie_b,Y)=-(2n+1)(\nabla_{\xi}\mathbb
B)(Y,\xi).
\end{gather}
The Ricci identities, equation \eqref{currr} in
Proposition~\ref{tech},  \eqref{cant2}, \eqref{tortrace} and the
symmetry of $L$ imply
\begin{multline}\label{rtt}
(\nabla^2_{\xi,e_b}L)(Ie_b,Y)-(\nabla^2_{e_b,\xi}L)(Ie_b,Y)=-
(\nabla_{Ie_b}\tau)(e_b,e_a)L(Y,e_a)
+\tau(e_b,Ie_a)(\nabla_{e_a}L)(e_b,Y)\\
+[(\nabla_{e_a}\tau)(e_b,Y)-(\nabla_Y)\tau)(e_b,e_a)]L(e_a,Ie_b).
\end{multline}
A small calculation taking into account \eqref{pp41},
\eqref{pp42}, \eqref{pp43}, \eqref{rtt} and using \eqref{tl}
yields
\begin{multline}\label{pp44}
(\nabla_{\xi}\mathbb B)(Y,\xi)-(\nabla_Y\mathbb
B)(\xi,\xi)-3[L(Y,Ie_a)+\tau(Y,e_a)]\mathbb
B(e_a,\xi)=-(trL)\mathbb B(IY,\xi)\\
+[(\nabla_{e_b}L(Y,Ie_a)+(\nabla_{e_b}\tau(Y,e_a)-(\nabla_{e_a}\tau(Y,e_b)+
(\nabla_Y)\tau(e_b,e_a)]L(Ie_b,e_a).
\end{multline}
Now, apply the already proven \eqref{inte} together with
\eqref{tl} to \eqref{pp44} to get the proof of \eqref{intehxi312}.
\end{proof}
Thus, the proof of Theorem~\ref{main2} i) is complete.
\subsection{The three dimensional case}
If the dimension is equal to 3 then it is easy to check that
$W^{pc}=0$ and the integrability conditions \eqref{inte} and
\eqref{inte1}  are trivially satisfied. Thus, the existence of a
smooth solution depends only on  the validity of \eqref{intexih11}
since the proof of Lemma~\ref{integr4} shows that
\eqref{intehxi312} follows from \eqref{intexih11} also in
dimension three. The next Lemma~\ref{thre} implies
Theorem~\ref{main2} ii).
\begin{lemma}\label{thre} If $n=1$ and
$F=0$  then \eqref{intexih11} holds.
\end{lemma}
\begin{proof}
Suppose $n=1$. Then $r(X,Y)=\rho(X,IY)=\frac{Scal}2 g(X,Y)$ and
\eqref{lll} yields
\begin{equation}\label{3l}
L(X,Y)=\frac{Scal}{16}g(X,Y)-\tau(X,IY).
\end{equation}
Apply \eqref{3l} to \eqref{pp6} to get
\begin{equation}\label{3xi}
2(\nabla_{\xi}L)(X,Y)=-(\nabla^2_{e_aIe_a}L)(X,Y) -
Scal.\tau(X,Y).
\end{equation}
The skew symmetric part of \eqref{intexih11} is satisfied because
$n=1$. Now, \eqref{pp2}, \eqref{3xi} and \eqref{3l} give that the
symmetric part of \eqref{intexih11} is equivalent to $F(X,Y)=0$.
\end{proof}
The proof  of Theorem~\ref{main2} is completed.


\section{A remark on the Cartan-Chern-Moser  theorem in the CR case}

A CR manifold is a smooth manifold $M$ of real dimension 2n+1,
with a fixed n-dimensional complex subbbundle $H$ of the
complexified tangent bundle $\mathbb CTM$ satisfying $ H \cap
\overline{H}=0$ and $[H,H]\subset H$. If we let $\mathbb H=Re\,
H\oplus\overline{H}$, the real subbundle $\mathbb H$ is equipped
with a formally integrable almost complex structure $J$. We assume
that $M$ is oriented and there exists a globally defined contact
form $\theta$ such that $\mathbb H=Ker\,\theta$. Recall that a
1-form $\theta$ is a contact form if the hermitian bilinear form
$2g(X,Y)=-d\theta(JX,Y)$ is non-degenerate. The vector field
$\zeta$ dual to $\theta$ with respect to $g$ and satisfying
$d\theta(\zeta, .)=0$ is called the Reeb vector field. A CR
manifold $(M,\theta,g)$ with fixed contact form $\theta$ is called
\emph{a pseudohermitian manifold}. In this case the 2-form
$d\theta_{|_{\mathbb H}}:=2\Omega$ is called the fundamental form.
Note that the contact form is determined up to a conformal factor,
i.e. $\bar\theta=\nu\theta$ for a positive smooth function $\nu$,
defines another pseudohermitian structure called pseudo-conformal
to the original one. A basic geometric tool in investigating
pseudohermitian structures is the Webster connection  $\Cr$
\cite{W,W1} (see also Tanaka \cite{T}).

The  Cartan-Chern-Moser results \cite{Car,ChM} are that the vanishing of
the pseudoconformal  invariant Chern-Moser tensor $S$ (resp. Cartan
invariant for $n=1$) is a necessary and sufficient condition a
non-degenerate  CR-hypersurface in $\mathbb{C}^{n+1}$, to be locally CR
equivalent to a hyperquadric in $\mathbb{C}^{n+1}$. A proof of these
results can be achieved working similarly to our proof of
Theorem~\ref{main2}. We outline below the crucial steps.

It is well known that a pseudohermitian manifold with flat Webster
connection (and zero torsion if $n=1$) is locally isomorphic to a
Heisenberg group. On the other hand, the Cayley transform is a
pseudo-conformal equivalence between the Heisenberg group with its
flat pseudo-hermitian structure and a hypersphere $g
_{\alpha\bar\beta}Z^\alpha Z^{\bar\beta}+W\bar W=1$ in
$\mathbb{C}^{n+1}$ \cite[p.223]{ChM}. It remains to show that the
vanishing of the Chern-Moser tensor, $S=0$, is a sufficient
condition a given pseudohermitian manifold to be locally
pseudoconformally flat provided the dimension is bigger than
three. In dimension three $S$ vanishes identically and the
sufficient condition remains only \eqref{crintexih11} below. The
scheme is formally very similar to that used in the proof of
Theorem~\ref{main2}. Namely in all formulas in the proof of
Theorem~\ref{main2} one formally replaces $I$ with $\sqrt{-1}J$
and $\xi$ by $\sqrt{-1}\zeta$.  We indicate below the most
important steps.

{ The superscript $cr$ means that the objects are taken with
respect to the Webster connection. In particular, the
pseudohermitian  Ricci 2-form $\rho^{cr}$ and the pseudohermitian
scalar curvature $Scal^{cr}$ are defined by $\quad
2\rho^{cr}(A,B)=g(\Rcr(A,B)\epsilon_a,J\epsilon_a), \quad
Scal^{cr}=r(\epsilon_a,\epsilon_a), \quad A,B\in \Gamma(TM).
\quad$ The Chern-Moser tensor $S$ can be obtained  from
\eqref{pinv} formally replacing $I$, $\rho$, $Scal$ and $\omega$
with $\sqrt{-1}J$, $\rho^{cr}$, $Scal^{cr}$ and  $\Omega$,
respectively.} The system of PDE which guaranties the flatness of
the pseudoconformal Webster connection and has to be solved is:

\begin{gather}\label{crsist1}
\Cr
dv(X,Y)=\hspace{11cm}\\\nonumber-C(X,Y)-dv(X)dv(Y)+dv(JX)dv(JY)+
\frac12g(X,Y)|\Cr v|^2 - dv(\zeta)\Omega(X,Y), \\\label{cradd1}
\Cr du(X,\zeta)=-\mathbb D(X,\zeta)-C(X,J\nabla
v)+\frac12dv(JX)|\Cr v|^2
-dv(X)du(\zeta),\hspace{2.5cm}\\\label{cradd2} \Cr
dv(\zeta,\zeta)=-\mathbb D(\zeta,\zeta) -\mathbb D(J \nabla
v,\zeta)+\frac14|\Cr v|^4-(dv(\zeta))^2,\hspace{4.5cm}
\end{gather}
where  the symmetric tensor $C(X,Y)$, $\mathbb D(X,\zeta)$ and $\mathbb
D(\zeta,\zeta)$ do not depend on the function $u$ and are determined by
\begin{equation}\label{crlll}
C(X,Y)=
-\frac1{2(n+2)}\rho^{cr}(X,JY)-\frac{Scal^{cr}}{8(n+1)(n+2)}g(X,Y)
+A(JX,Y)
\end{equation}
\begin{gather}\label{crbes}
\begin{aligned}
(\Cr_{e_a}C)(Je_a,JX)=-(2n+1)\mathbb D(JX,\zeta),
\end{aligned}\\\label{crbst}
 \mathbb D(\zeta,\zeta)=-\frac{1}{2n}[(\Cr_{e_a}\mathbb
D)(Je_a,\zeta)+C(e_b,Je_a)C(Je_b,e_a)],
\end{gather}
and the symmetric tensor $A(X,Y)$ is the pseudohermitian torsion
\cite{W,W1,L2}. }

The integrability conditions for the overdetermined system
\eqref{crsist1}-\eqref{cradd2} are:
\begin{gather}\label{crinte}
(\Cr_ZC)(X,Y)-(\Cr_XC)(Z,Y)=-\Omega(Z,Y)\mathbb D(X,\zeta)+
\Omega(X,Y)\mathbb D(Z,\zeta)-2\Omega(Z,X)\mathbb D(Y,\zeta);
\\\label{crinte1}
(\Cr_Z\mathbb D)(X,\zeta)-(\Cr_X\mathbb
D)(Z,\zeta)+C(Z,JC(X))-C(X,JL(Z))= -2\mathbb
D(\zeta,\zeta)\Omega(Z,X);\\
\label{crintexih11} (\Cr_X\mathbb
D)(Y,\zeta)-(\Cr_{\zeta}C)(X,Y)=\hspace{8cm}\\\nonumber
\hspace{4cm} C(Y,JC(X))+A(X,C(Y))+A(Y,C(X))-\mathbb
D(\zeta,\zeta)\Omega(X,Y);\\
\label{crintehxi312} (\Cr_{\zeta}\mathbb D)(X,\zeta)-(\Cr_X\mathbb
D)(\zeta,\zeta)-2\mathbb D(e_a,\zeta)C(X,Je_a)+ A(X,e_a)\mathbb
D(e_a,\zeta)=0.
\end{gather}

As in the proof of Theorem~\ref{main2} i) we can see that the vanishing of
the Chern-Moser tensor, $S=0$, implies the validity of the integrability
conditions \eqref{crinte}-\eqref{crintehxi312} provided $n>1$.

For $n=1$ the Chern-Moser tensor is always zero and, following the
proof of Theorem~\ref{main2} ii),  one checks  that the
integrability conditions \eqref{crinte} and \eqref{crinte1} are
trivially satisfied and \eqref{crintehxi312} is  a consequence of
\eqref{crintexih11}. To have a smooth solution to the system
\eqref{crsist1}-\eqref{cradd2} one has to have \eqref{crintexih11}
which is equivalent to the vanishing of the symmetric (0,2) tensor
$F^{car}$  defined on $\mathbb H$ by
\begin{gather}\label{crtreeF}
F^{car}(X,Y)=(\Cr d(Scal^{cr}))(X,JY)+(\Cr
d(Scal^{cr}))(Y,JX)\hspace{3cm}\\\nonumber \hspace{3.5cm}+16((\Cr)
^2_{Xe_a}A)(Y,e_a) +16((\Cr)^2_{Ye_a}A)(X,e_a)
+48((\Cr)^2_{e_aJe_a}A)(X,JY)\\\nonumber+36Scal^{cr}A(X,Y)+3g(X,Y)(\Cr
d(Scal^{cr}))(e_a,Je_a).
\end{gather}
{ Let us remark that the vanishing of $F^{car}$ is equivalent to
the vanishing of the Cartan curvature, cf.  \cite[Theorem
12.3]{Tann}.}
\begin{cor}
A 3-dimensional Sasakian manifold $(M,\theta,g,\zeta)$ is locally
pseudoconformally equivalent to the three dimensional Heisenberg
group if and only if its Riemannian scalar curvature $Scal^g$
satisfies
\begin{gather}\label{scalg}
(\nabla^g d(Scal^g))(X,JY)+(\nabla^g d(Scal^g))(Y,JX)=0,
\end{gather}
{ where $\nabla^g$ is the Levi-Civita connection of $g$.}
\end{cor}
\begin{proof}
It is well known that a pseudohermitian structure is Sasakian,
i.e. its Riemannian cone is K\"ahler, exactly when the Webster
torsion vanishes, $A=0$. { In particular, the Bianchi identities
imply $\zeta(Scal^{cr})=0$.} { Then the second and the third lines
in \eqref{crtreeF} disappeared in view of \eqref{symdh}. On the
other hand, for a Sasaki manifold, we have
$\Cr_XY=\nabla^g_XY+g(JX,Y)\zeta$ and the Riemannian scalar
curvature and the scalar curvature of the Webster connection
differ by an additive constant depending on the dimension,
$2Scal^{cr}=Scal^g+2n$ (see e.g. \cite{DT}). Now, \eqref{scalg}
becomes equivalent to \eqref{crtreeF}. Hence, $(M,\theta)$ is
locally pseudoconformally flat.}
\end{proof}

\section{The ultrahyperbolic Yamabe equation}\label{Yamabe}

Recall that the CR Yamabe problem is to determine if there exists
a pseudohermitian structure compatible with a given CR structure
such that the pseudohermitian scalar, i.e. the scalar curvature of
the Webster connection is  constant. If the CR structure is
strongly pseudo-convex, i.e. the Levi form is negative definite
then the CR Yamabe problem reduces to a subelliptic PDE which can
be solved  on a compact manifold \cite{JL1}.

Similarly to the CR case one can pose a Yamabe type problem for a
para CR manifold. Namely, given a para CR structure is there a compatible 
para hermitian structure such that the
scalar curvature of the canonical connection is a constant.

In the case when the Levi form of a given CR structure  has
neutral signature of type (n,n) then the CR Yamabe equation is of
the same type as the para CR Yamabe equation \cite{L2}, i.e. one
has to consider the sub ultrahyperbolic equation \eqref{e:conf
change scalar curv} with respect to the Webster connection where
$\bar s=const.$.

Here we show an explicit formula for the regular
part of a solution to the ultrahyperbolic Yamabe equation on $\GP$
\begin{equation}\label{e:Yamabe}
\lap \varphi\ \equiv\ \sum_{k = 1}^n \left ( U^2_k \ - \ V^2_k\right
) \varphi \ =\ -\varphi^{2^*-1},
\end{equation}
where $2^*-1=(Q+2)/(Q-2)=(n+2)/n$ with $Q=2n+2$ the homogenous
dimension of the group. When  $n=2m$ \eqref{e:Yamabe} coincides with
the Yamabe equation on the Heisenberg group of (real) signature
$(2m,2m)$ defined by the quadric
\[
Q=\{ (z,w)\in\mathbb{C}^n\times\mathbb{C}:\ \text{Im}\, w=H(z,z)
 \},
\]
where $H(z,z)=\sum_{j=1}^m \bigl(\,  z_j
\overline{z}_j'-z_{j+m}\overline{z}_{j+m}'\, \bigr )$, with the
natural group structure
\[
( z'',w'')\ =\ (z',w')\circ(z, w)\ =\ \bigl(z'\ +\ z, w' +\ w\
+2\text{Im}\, H(z',z)\bigr ).
\]
The left-invariant horizontal vector fields are given by
\begin{align*}
& X_j=\dxj-2y_j\ddt, \quad Y_j=\dyj + 2x_j\ddt\\
& X_{m+j}=\dxmj+2y_{j+m}\ddt, \quad Y_{m+j}=\dymj-2x_{j+m}\ddt,
\quad j=1,\dots,m,
\end{align*}
while the left invariant contact form with corresponding metric, for
which the above vector fields are an orthonormal frame, is given by
\begin{equation*}
\theta  = \frac 12\,dt + \Sigma_{j=1}^m (y_j\,dx_j-x_j\,dy_j) -
\Sigma_{j=1}^m (y_{j+m}\,dx_{j+m}-x_{j+m}\,dy_{j+m})
\end{equation*}
so that $
g(X_j,X_j)=-g(Y_j,Y_j)=-g(X_{j+m},X_{j+m})=g(Y_{j+m},Y_{j+m})=1,\quad
j=1,\dots,m. $

By the Cartan-Chern-Moser result, the above quadric is the flat CR
structure of (hermitian) signature $(m,m)$. Henceforth, for
$\vu\in\mathbb{R}^n$, $\vu=(u_1,\dots,u_n)$ we set $|\vu|=\bigl
(u_1^2+\dots u_n^2 \bigr )^{1/2}$. We observe 
\begin{prop}\label{T:Yamabe}
Let $\GP$ be the  Heisenberg group of topological dimension $2n+1$.
For every $\epsilon> 0$ the function
\begin{equation}\label{Ksigma}
\varphi_\epsilon(\vu,\vv,t)\ =\ \left(\frac{4n^2 \epsilon^2}{
(\epsilon^2 + \abs{\vu}^2-\abs{\vv}^2)^2 - t^2}\right)^{
\frac{n}{2}},\quad\quad\quad g\in \bG,
\end{equation}
is a solution of the ultrahyperbolic Yamabe equation
\eqref{e:Yamabe} on the set where $\abs{\epsilon^2 +
\abs{\vu}^2-\abs{\vv}^2} \not= \abs{t}$.
\end{prop}
\begin{proof}
Let $f=\left( (1 + \abs{\vu}^2-\abs{\vv}^2)^2 - t^2\right )^{
-\frac{n}{2}}$. After a straightforward calculation we find $\lap f
= -4n^2 f^{2^*-1}$, which implies easily the equation for
$\varphi_1$. Furthermore, using the dilations on the group
$\delta_\lambda (\vu,\vv,t)= ( \lambda u, \lambda v,\lambda^2 t)$
we have that the function $f_\lambda(\vu,\vv,t)= \lambda ^{n/2}f(
\lambda u, \lambda v,\lambda^2 t)$ satisfies the same equation as
$f$, which implies the equation for $\varphi_\epsilon$ by taking
$\epsilon=1/\lambda$.
\end{proof}
Since the ultra-hyperbolic Yamabe equation is invariant under
translations it follows that we can construct other solutions,  each
being a regular function on a corresponding set.  The question
whether there is a global solution, in the sense of distributions,
will not be considered here. In this respect we note that
\cite{Tie}, \cite{MR} found the fundamental solution of the
ultra-hyperbolic operator in the left-hand side of \eqref{e:Yamabe}.

It should be pointed out that there is a correspondence between
the regular part of solutions to partial differential equations on
the hyperbolic Heisenberg group and solutions of partial
differential equations on the Heisenberg group. Let  $X_{k}=${\
}$\frac{\partial }{\partial x_{k}}-2y_{k} \frac{\partial
}{\partial s},$ $Y_{k}=\frac{\partial }{\partial y_{k}}+2x_{k}
\frac{\partial }{\partial s}$ be the horizontal left invariant
vector fields on the standard Heisenberg group. Note that the
difference between this group and the hyperbolic Heisenberg group
is in the metric, while the groups are identical. Given a function
$f(x,y,t)$, $t\in\mathbb{R}, \ x,y\in\mathbb{R}^n$, let
$g(u,v,t)=f(it,u,iv)$, which  could be  a complex valued function
even when $f$ is real-valued. Since
\begin{equation}\label{e:formal change}
(X_kf) (\vx,\vy,s)=(U_kg) (\vu,\vv,t),\qquad (Y_kf)
(\vx,\vy,s)=-i(V_kg) (\vu,\vv,t)
\end{equation}
we have $\sum_{k = 1}^n \left ( X^2_k \ + \ Y^2_k\right ) f =\sum_{k
= 1}^n \left ( U^2_k \ - \ V^2_k\right ) f $. In particular,
solutions of the Yamabe equation on the Heisenberg group turn into
 solutions of the ultra-hyperbolic Yamabe equation on the hyperbolic  Heisenberg
group outside a corresponding singular set.

Consider the (standard) Heisenberg group of dimension $2n+1$ with
typical point $(z,s)$, $z\in\Cn,\, s\in\mathbb{R}$, and let
$A=\nz^2+it$. The inversion of the point $(z,s)$ is given by
\begin{equation}
(z',s')\overset {def}{=}\bigl( -\frac {z}{\Ab}, -\frac {s}{A\Ab}
\bigr),
\end{equation}
which can also be written  in real coordinates as
\[
 \vx'=\vx'(\vx,\vy,s), \qquad \vy'=\vy'(\vx,\vy,s), \qquad s'=s'(\vx,\vy,s).
\]
If $B\overset {def}{=}\abs{z'}^2+is'$, then $AB=1$ as $B=\frac
{\nz^2}{A\Ab}-\frac {is}{A\Ab}=\frac {\Ab}{A\Ab}=\frac {1}{A}.$

 Based on the above mentioned formal substitution and the preceding
 paragraph we define an inversion on the hyperbolic Heisenberg group  as follows.
 Let $\Xi=\{ p=(\vu,\vv,t)\in\GP\, :\
\left\vert \abs{\vu}^2-\abs{\vv}^2\right\vert = \abs{t}\}$ and
$p=(t,u,v)\in\GP\setminus\Xi$. We define the inversion on the
hyperbolic Heisenberg group letting
\begin{equation}\label{e:real inversion}
\vu'=\vx'(\vu,i\vv,it), \qquad \vv'=-i\vy'(\vu,i\vv,it), \qquad
t'=-is'(\vu,i\vv,it)
\end{equation}
using the real form of the "standard" inversion on the Heisenberg
group. In other words, for $k=1,\dots,n$ we have
\begin{equation}\label{e:inversion}
u_k'=-\frac {{\nw}u_k+tv_k}{{\nw}^2-t^2}, \quad v_k'=-\frac
{{\nw}v_k+tu_k}{{\nw}^2-t^2}, \quad t'=-\frac {t}{{\nw}^2-t^2},
\end{equation}
which defines a point $p'=(\vu',\vv',t')\in\GP\setminus\Xi$, or,
using $\vw=\vu+e\vv$ with $e^2=1$, $|\vw|^2=|\vu|^2-|\vv|^2$,
\[
\vw'\equiv \vu'+e\vv'=-\frac {\vw}{|\vw|^2-et}\qquad t'=-\frac
{t}{|\vw|^4-t^2}.
\] This map will be called the inversion of $\GP$ centered at $\Xi$.
The inverse transformation is found by taking into account that the
inversion is an involution.

Recall, see\cite{K}, that for a function $f(z,t)$ defined on a
domain $\Om$ in the Heisenberg group  we define the Kelvin transform
${f\st}$ on the image ${\Om\st}$ of $\Om$ under the inversion by the
following formula
\begin{equation*}
{f\st}\overset {def}{=} A^{\frac {n}{2}}\Ab^{\frac {n}{2}}f \text {
i.e. } {f\st}\nB^n=f.
\end{equation*}
Thus, using the preceding considerations we can define a Kelvin
transform on the hyperbolic Heisenberg group as follows
\begin{equation*}
(\mathcal{K}\varphi) (\vu,\vv,t)\\
=\bigl( (|\vu|^2-|\vv|^2)^2 - t^2 \bigr)^{-n/2}\,
\varphi(\vu',\vv',t'),
\end{equation*}
where $(\vu'\vv',t')$  are given by \eqref{e:inversion}.  Given a
function $\varphi(u,v,t)$ we consider
$\psi(x,y,s)=\varphi(x,-iy,-is)$. Thus, using $s=it$, $\vx=\vu$ and
$\vy=i\vv$, we have from  \eqref{e:real inversion}
\begin{multline*}
\varphi(\vu',\vv',t')=\psi(\vu',i\vv',it')=\psi\bigl
(\vx'(\vu,i\vv,it),\vy'(\vu,i\vv,it),s'(\vu,i\vv,it) \bigr)\\=\psi
\bigl (\vx'(\vx,\vy,s),\vy'(\vx,\vy,s),s'(\vx,\vy,s) \bigr),
\end{multline*}
which shows that the (hyperbolic) Kelvin transform of $\varphi$
corresponds to the ("standard" Heisenberg ) Kelvin transform of
$\psi$. Due to \eqref{e:formal change} and the properties of the
Kelvin transform on the Heisenberg group (in fact any group of
Iwasawa type), cf. \cite{CDKR} and \cite{GV}, the hyperbolic Kelvin
transform preserves the ultra-hyperbolic  functions, i.e., solutions
of
\[ \lap \varphi\ \equiv\ \sum_{k = 1}^n \left ( U^2_k \ - \
V^2_k\right ) \varphi \ =0 \] and the solutions of the
ultra-hyperbolic Yamabe equation. Note that the Kelvin transform of
the functions in \eqref{Ksigma} is given by the same formula.

\end{document}